\documentclass{elsarticle}  %
\usepackage{graphicx}
\usepackage{amsmath}
\usepackage{amssymb}
\usepackage{enumerate}
\usepackage{theorem}
\usepackage{rawfonts}
\usepackage{setspace}

\usepackage{arydshln}
\usepackage{lscape}

\usepackage{float}
\usepackage{color}
\usepackage{multicol}

\usepackage[all]{xy}
\xyoption{poly} \xyoption{graph} \xyoption{arc} \xyoption{web}
\xyoption{2cell}

\begingroup
\newtheorem{thrm}{Theorem}[section]
\newtheorem{crry}[thrm]{Corollary}

\newtheorem{lemm}[thrm]{Lemma\hskip 1mm}

\newtheorem{dfn} [thrm]{Definition}

\endgroup

\newenvironment{pf}{\noindent\textbf{Proof}}{\hspace*{\fill}$\square$\\[6pt]}

\begin{document}

\title{On connected simple graphs and their degree sequences}

\author{Jonathan McLaughlin}
\address{Department of Mathematics, St. Patrick's College, Dublin City University, Dublin 9, Ireland }
\ead{jonny$\_\:$mclaughlin@hotmail.com}

\parindent=0cm


\begin{abstract} This note describes necessary and sufficient conditions for a sequence of positive integers to be the degree sequence of a connected simple graph. Conditions are also given under which a sequence is necessarily connected i.e. the sequence can only be realised as a connected graph. A matrix is introduced whose non-empty entries  partition the set of connected graphs. The note concludes with a result relating the number of edges in a simple graph to the connectedness of the graph. 
\end{abstract}


\maketitle

\parindent=0cm

\section{Introduction}
A finite sequence of non-negative integers is called a connected sequence if it is the degree sequence of some finite simple connected  graph. Necessary and sufficient conditions for a sequence of non-negative integers to be connected are implicit in Hakimi \cite{Hk} and these conditions have been stated explicitly by the author in \cite{Me15}. The current note takes a constructive approach to connected degree sequences. An iterative construction for connected graphs, where each intermediate graph is connected, is described before the set of connected graphs is partitioned into equivalence classes each of which is placed separately in a unique entry of a matrix $P_{\mathcal{G}}$. Necessary and sufficient conditions are then given for a sequence to be connected and/or necessarily connected. The note concludes with a result relating the connectedness of a simple graph to the number of edges in the graph.

\section{Preliminaries }

Let $G=(V_{G},E_{G})$ be a graph where $V_{G}$ denotes the vertex set of $G$ and $E_{G}\subseteq [V_{G}]^{2}$ denotes the edge set of $G$ (and $[V_G]^2$ is the set of all $2$-element subsets of $V_G$).  An edge $\{a,b\}$ is denoted $ab$. A graph is finite when $|V_{G}|<\infty$ and $|E_{G}|<\infty$, where $|X|$ denotes the cardinality of the set $X$. The union of graphs $G$ and $H$ i.e. $(V_{G}\cup V_{H}, E_{G}\cup E_{H})$, is denoted $G\cup H$. By a slight abuse of notation, $ab\cup G$ is understood to be the graph $(\{a,b\},\{ab\})\cup (V_G, E_G)$. A graph is simple if it contains no loops (i.e. $aa\not\in E_{G}$) or parallel/multiple edges (i.e. $\{ab,ab\}\not\subseteq E_{G}$). The {\it degree} of a vertex $v$ in a graph $G$, denoted $deg(v)$, is the number of edges in $G$ which contain $v$. A {\it path} is a graph with $n$ vertices in which two vertices, known as the {\it endpoints}, have degree $1$ and $n-2$ vertices have degree $2$. A graph is {\it connected} if there exists at least one path between every pair of vertices in the graph. A {\it tree} is a connected graph with $n$ vertices and $n-1$ edges. $K_{n}$ denotes the {\it complete graph} on $n$ vertices. All basic graph theory definitions can be found in standard texts such as \cite{BM}, \cite{D} or \cite{GG}. All graphs in this work are assumed to be simple, undirected and finite. \\

A finite sequence $s=\{s_1,...,s_n\}$ of non-negative integers is called {\it graphic} if there exists a finite simple graph with vertex set  $\{v_1,..., v_n\}$ such that $v_i$ has degree $s_i$ for all $i=1,...,n$. Necessary and sufficient conditions for a sequence of non-negative integers to be graphic were first described by Erd\"os and Gallai in \cite{EG} and are extended in works by Choudum \cite{Ch}, Sierksma $\And$ Hoogeveen \cite{SH} and Tripathi et al. \cite{TVW10}, \cite{TV03}, \cite{TV07}. In this note all sequences have {\it positive} terms as the only connected graph which has a degree sequence containing a zero is $(\{v\},\{\})$. The maximum degree of a vertex in $G$ is denoted $\Delta_G$ and the minimum degree of a vertex in $G$ is denoted $\delta_G$. Given a graph $G$ then the {\it degree sequence} $d(G)$ is the monotonic non-increasing sequence of degrees of the vertices in $V_G$. This means that every graphical sequence $s$ is equal to the degree sequence $d(G)$ of some graph $G$ (subject to possible rearrangement of the terms in $s$).

\begin{dfn} A finite sequence $s=\{s_1,...,s_n\}$ of positive integers is called {\it connected} if there exists a finite simple connected graph with vertex set  $\{v_1,..., v_n\}$ such that $deg(v_i)= s_i$ for all $i=1,...,n$.
\end{dfn}

\section{Constructing connected graphs}
\begin{dfn}\label{hpathdfn} Given a graph $H=(V_H,E_H)$, then $ab=(\{a,b\},\{ab\})$ is called a {\it closed} $H${\it-edge} if $ab\cap H=(\{a,b\},\{\})$ or an {\it open} $H${\it-edge} if $ab\cap H=(\{a\}, \{\})$ (or equivalently, if $ab\cap H=(\{b\}, \{\})$). \end{dfn}

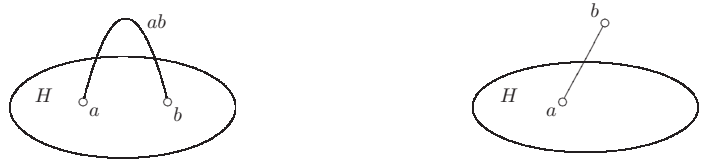
\begin{figure}[h]
\begin{center}
\scalebox{0.75}{$\begin{xy}\POS (0,5) *\cir<2pt>{} ="a" *+!UL{a},
(15,5) *\cir<2pt>{} ="b" *+!UL{b},
  (-7,5)*+!{H},
  (13,18)*+!{ab},

\POS "a" \ar@/^3.5pc/@{-}  "b",

\POS(7,4),  {\ellipse(20,9)<>{}},

\POS 
(85,5) *\cir<2pt>{} ="a1" *+!UR{a},
(92.5,19) *\cir<2pt>{} ="u" *+!DR{b},
  (75.5,5)*+!{H},

\POS "a1" \ar@{-}  "u",

\POS(89,4),  {\ellipse(20,8)<>{}},

 \end{xy}$}

\caption{A closed $H$-edge and an open $H$-edge, respectively. }
\label{Cases2}
\end{center}
\end{figure}

Lemma \ref{connLemm} now describes an iterated construction of a connected graph, beginning with $K_1=(\{v\},\{\})$, where every intermediate graph is also connected. 

\begin{lemm}\label{connLemm} A graph $G$ is connected if and only if $G$ can be constructed from $K_1$ by successively adding open or closed $H$-edges to subgraphs $H$ which have already been constructed.
\end{lemm}

\begin{pf}($\Rightarrow$) Given a connected graph $H$, the smallest of which is $K_1$, then $H\cup ab$, where $ab$ is either an open or closed $H$-edge, is also connected. Every connected graph $G$ can be constructed in two steps. In the first step, $T$, a spanning tree of $G$, is built from $K_1$ by a series of unions involving  open $H$-edges. The second step adds the edges in $E_{G}\setminus E_T$ and is achieved using a series of unions involving closed $H$-edges. \\ 

($\Leftarrow$) Assume that $G$ is connected, then $G$ contains a graph which is isomorphic to $K_{1}$ and hence also contains a maximal subgraph $H$ which is constructible as per the statement of the result. Suppose that $G\neq H$ and consider the subgraph $G\setminus H$ contained in $G$. Observe that $G\setminus H$ cannot contain an edge $ab$ where $ab \not\in E_{H}$ and $a,b\in V_{G}$ as such an edge would be a closed $H$-edge and this would contradict the maximality of $H$. So, $G\setminus H$ must contain at least one vertex $u$ such that $u\not\in V_{H}$. As $G$ is connected then it is possible to choose a vertex $u\in V_{G\setminus H}$ such that there is at least one edge $ua$ in $G$ where $a\in V_{H}$ . However, $ua$ would then be an open $H$-edge which again contradicts the maximality of $H$. Hence $H=G$. 
 \end{pf}

\section{Listing all connected graphs}\label{s4}

Let $\mathcal{G}$ be the set of all (unlabelled) connected simple graphs. Define a graph equivalence relation $\sim$ as $G\sim H$ if and only if $|V_G|=|V_H|$ and $|E_G|=|E_H|$.\\

Consider a matrix $P=\big( p_{i,j}\big)_{i,j\in \mathbb{N}_0}$. Each non-empty entry $p_{i,j}\in P$ contains a single equivalence class (with respect to $\sim$) with the entry $p_{0,0}$ containing $K_1$. \\

The equivalence classes are arranged in $P$ as follows:
\begin{itemize}
\item Each $G\in p_{i,j}$ is a maximal proper subgraph of some $G'\in p_{i,j+1}$ such that $V_{G}\subset V_{G'}$. 
\item Each $G\in p_{i,j}$ is a maximal proper subgraph of some $G''\in p_{i+1,j}$ such $V_{G}= V_{G''}$ and $E_{G}\subset E_{G''}$. 
\item If $p_{i,j}$ contains a maximal complete graph, then $p_{i+1,j},\; p_{i+2,j}, \; ...$ are empty.
\end{itemize}

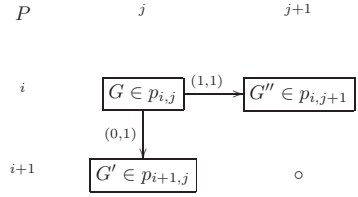
\begin{figure}[h]
\begin{center}\scalebox{0.75}{$
\xymatrix{P & ^{j} & ^{j+1} &    \\
 ^{i}   &    *+[F]{G\in p_{i,j} }  \ar@{->}^<<<<<{(1,1)}[r] \ar@{->}_{(0,1)}[d] & *+[F]{G''\in p_{i,j+1}} &    \\
^{i+1}   &     *+[F]{G'\in p_{i+1,j} } &  \circ &   \\ 
}$ }  
\caption{The relationship between graphs in adjacent entries in $P$. }
\label{Poly2}
\end{center}
\end{figure}

Such an arrangement of equivalence classes within $P$ encodes the construction described in Lemma \ref{connLemm} in the following manner: 
\begin{itemize}
\item Assuming that $p_{i,j}$ does not contain a maximal complete graph, then the union of some $G\in p_{i,j}$ and a closed $H$-edge results in some $G'\in p_{i+1,j}$ where $G\subset G'$ with $|V_{G}|= |V_{G'}|$ and $|E_{G}|+1= |E_{G'}|$. This is signified by the $(0,1)$ label in Figure \ref{Poly2}.
\item The union of some $G\in p_{i,j}$ and an open $H$-edge results in some $G''\in p_{i,j+1}$ where $G\subset G''$ with $|V_{G}|+1= |V_{G''}|$ and $|E_{G}|+1= |E_{G''}|$. This is signified by the $(1,1)$ label in Figure \ref{Poly2}. 
\end{itemize}

For each $G\in p_{i,j}$ it is possible to state $|V_G|$ and $|E_G|$ in terms of $i$ and $j$. As $p_{0,0}$ contains $K_1$ and each column of $P$ contains graphs with the same number of vertices, then column $j$ contains graphs with $j+1$ vertices. \\

The uppermost row of $P$ (i.e. row $0$) contains minimally connected graphs i.e. trees, and so graphs in entry $p_{0,j}$ have $j+1$ vertices and $|V_G|-1=j$ edges. Observe that graphs in $p_{1,j}$ (if they exist) contain the same number of vertices as graph in $p_{0,j}$ but an additional edge, hence, graphs in $p_{1,j}$ have $j+1$ edges. It follows that graphs in $p_{i,j}$ (if they exist) have $j+1$ vertices and $j+i$ edges.\\

It is now possible to define $P_{\mathcal{G}}$, the partition matrix of $\mathcal{G}$.

\begin{dfn} Let $\mathcal{G}$ be the set of all connected simple graphs then the partition matrix of $\mathcal{G}$ is the matrix \[P_{\mathcal{G}} :=\big( p_{i,j}\big)_{i,j\in \mathbb{N}_0}\]
such that \[ p_{i,j} :=\{G\in \mathcal{G} \; \mid \;  |V_{G}|=j+1\;\; {\textrm and } \;\; |E_{G}|= i +j  \}.\]
\end{dfn}

A portion of $P_{\mathcal{G}}$ in shown in Figure \ref{Gamma}.

\begin{figure}[h]
\[ \scalebox{0.95}{$P_{\mathcal{G}}= \left( \begin{array}{ccccc}
p_{0,0} & p_{0,1} & p_{0,2} & p_{0,3} & \cdots \\
p_{1,0} & p_{1,1} & p_{1,2} & p_{1,3} & \cdots \\
p_{2,0} & p_{2,1} & p_{2,2} & p_{2,3} & \cdots \\
p_{3,0} & p_{3,1} & p_{3,2} & p_{3,3} & \cdots \\
p_{4,0} & p_{4,1} & p_{4,2} & p_{4,3} & \cdots \\
\vdots & \vdots & \vdots &\vdots & \ddots \\
\end{array} \right) =  \left( \begin{array}{ccccc}
  $\begin{xy} \POS  (0,0) *\cir<2pt>{} ="a" *+!UR{}, \end{xy}$ & $\begin{xy} \POS  (0,0) *\cir<2pt>{} ="a" *+!UR{}, (5,0) *\cir<2pt>{} ="b" *+!UR{}, \POS "a" \ar@{-}  "b", \end{xy}$  & $\begin{xy} \POS  (0,0) *\cir<2pt>{} ="a" *+!UR{}, (5,0) *\cir<2pt>{} ="b" *+!UR{},  (5,5) *\cir<2pt>{} ="c" *+!UR{}, \POS "a" \ar@{-}  "b", \POS "c" \ar@{-}  "b", \end{xy}$ & $\begin{xy} \POS  (0,0) *\cir<2pt>{} ="a" *+!UR{}, (5,0) *\cir<2pt>{} ="b" *+!UR{},  (5,5) *\cir<2pt>{} ="c" *+!UR{}, (0,5) *\cir<2pt>{} ="d" *+!UR{}, (10,0) *\cir<2pt>{} ="a1" *+!UR{}, (15,0) *\cir<2pt>{} ="b1" *+!UR{},  (15,5) *\cir<2pt>{} ="c1" *+!UR{}, (10,5) *\cir<2pt>{} ="d1" *+!UR{}, \POS "a" \ar@{-}  "b", \POS "b" \ar@{-}  "c",  \POS "c" \ar@{-}  "d", \POS "a1" \ar@{-}  "b1", \POS "b1" \ar@{-}  "c1", \POS "b1" \ar@{-}  "d1",  \end{xy}$  & \cdots  \\
\varnothing & \varnothing & $\begin{xy} \POS  (0,0) *\cir<2pt>{} ="a" *+!UR{}, (5,0) *\cir<2pt>{} ="b" *+!UR{},  (5,5) *\cir<2pt>{} ="c" *+!UR{}, \POS "a" \ar@{-}  "b", \POS "c" \ar@{-}  "b", \POS "a" \ar@{-}  "c", \end{xy}$  & $\begin{xy} \POS  (0,0) *\cir<2pt>{} ="a" *+!UR{}, (5,0) *\cir<2pt>{} ="b" *+!UR{},  (5,5) *\cir<2pt>{} ="c" *+!UR{}, (0,5) *\cir<2pt>{} ="d" *+!UR{}, (10,0) *\cir<2pt>{} ="a1" *+!UR{}, (15,0) *\cir<2pt>{} ="b1" *+!UR{},  (15,5) *\cir<2pt>{} ="c1" *+!UR{}, (10,5) *\cir<2pt>{} ="d1" *+!UR{}, \POS "a" \ar@{-}  "b", \POS "b" \ar@{-}  "c", \POS "c" \ar@{-}  "d", \POS "a" \ar@{-}  "d", \POS "a1" \ar@{-}  "b1", \POS "b1" \ar@{-}  "c1", \POS "b1" \ar@{-}  "d1",  \POS "d1" \ar@{-}  "c1", \end{xy}$  &  \cdots\\
\varnothing & \varnothing & \varnothing & $\begin{xy} \POS  (0,0) *\cir<2pt>{} ="a" *+!UR{}, (5,0) *\cir<2pt>{} ="b" *+!UR{},  (5,5) *\cir<2pt>{} ="c" *+!UR{}, (0,5) *\cir<2pt>{} ="d" *+!UR{},  \POS "a" \ar@{-}  "b", \POS "b" \ar@{-}  "c", \POS "c" \ar@{-}  "d", \POS "a" \ar@{-}  "d",  \POS "b" \ar@{-}  "d", \end{xy}$ &  \cdots \\
\varnothing & \varnothing & \varnothing & $\begin{xy} \POS  (0,0) *\cir<2pt>{} ="a" *+!UR{}, (5,0) *\cir<2pt>{} ="b" *+!UR{},  (5,5) *\cir<2pt>{} ="c" *+!UR{}, (0,5) *\cir<2pt>{} ="d" *+!UR{},  \POS "a" \ar@{-}  "b", \POS "b" \ar@{-}  "c", \POS "c" \ar@{-}  "d", \POS "a" \ar@{-}  "d",  \POS "b" \ar@{-}  "d",   \POS "a" \ar@{-}  "c", \end{xy}$  &  \cdots \\
\varnothing & \varnothing & \varnothing & \varnothing  &  \cdots \\
\vdots & \vdots & \vdots & \vdots &  \ddots \\  \end{array} \right) $} \] 
\caption{$P_{\mathcal{G}}$, the partition matrix of $\mathcal{G}$.}
\label{Gamma}
\end{figure}

\begin{lemm}\label{CC} Given a connected graph $G$ such that $|V_G|=v$ and $|E_G|=e$ then $G$ can be constructed from $K_1$ by successively adding $v-1$ open $H$-edges followed by $e-v+1$ closed $H$-edges to subgraphs $H$ which have already been constructed.
\end{lemm}
\begin{pf} Let $G$ be a connected graph with $|V_G|=v$ and $|E_G|=e$. Consider the partition matrix $P_{\mathcal{G}}$ along with Figure \ref{Poly2}. Recall that $p_{0,0}$ contains $K_1$. Each graph in $p_{0,j+1}$ is the union of a graph in $p_{0,j}$ and an open $H$-edge and so after $v-1$ of these unions the resultant graphs contain $v$ vertices and $v-1$ edges and one such graph is a spanning tree of $G$. Each graph in $p_{i+1,j}$ is the union of a graph in $p_{i,j}$ and a closed $H$-edge and so after $e-(v-1)$ of these unions the graph $G$ containing $v$ vertices and $e$ edges has been constructed. 
\end{pf}
Lemma \ref{CC} describes a {\it canonical} construction for any connected graph $G$ but it is worth noting that the order in which the open and closed $H$-edges are added is not arbitrary. Observe in Figure \ref{Gamma} that the $4$-cycle in $p_{1,3}$ must be constructed from $K_1$ using unions with three open $H$-edges followed by a union with a closed $H$-edge as a $4$-cycle is not the result of the union of any graph in $p_{1,2}$ and an open $H$-edge.

\section{Degree sequences and graph partitions}\label{s5}

Given a sequence of positive integers $s=\{s_1,...,s_n\}$ then define the {\it associated pair of $s$}, denoted $(\varphi(s),\epsilon(s))$, to be the pair $(n, \frac{1}{2}\sum\limits_{i=1}^{n}s_i)$. Where no ambiguity can arise,  $(\varphi(s),\epsilon(s))$ is simply denoted $(\varphi,\epsilon)$. 
 
\begin{lemm}\label{Ent} Given a sequence $s$ with associated pair $(\varphi,\epsilon)$ such that $s=d(G)$ for some connected graph $G\in \mathcal{G}$, then $G$ is contained in $p_{\;\epsilon-\varphi +1, \; \varphi -1}$ in $P_{\mathcal{G}}$. 
\end{lemm}

\begin{pf} As $s$ is the degree sequence $d(G)$ of some graph $G\in \mathcal{G}$ then $|V_G|$ is the number of terms in $s$ which is $\varphi$ and $|E_G|$ is half the sum of the degrees of all vertices in $G$ which is exactly $\epsilon$. Recall that each entry $p_{i,j}$ in $P_{\mathcal{G}}$ is defined as $\{G\in \mathcal{G} \; \mid \;  |V_{G}|=j+1\;\; {\textrm and } \;\; |E_{G}|= i+j  \} $. By rearranging $|V_{G}|=\varphi=j+1$ and substituting into $|E_{G}|=\epsilon= i+j $ then $j=\varphi-1$ and $i=\epsilon -\varphi +1$ hence $s=d(G)$ for some $G \in \mathcal{G}$ is contained in $p_{\;\epsilon -\varphi +1,\;\varphi-1}$ in $P_{\mathcal{G}}$.  
\end{pf}

\begin{lemm}\label{EmpEnt} The non-empty entries in column $j=\varphi -1$ of $P_{\mathcal{G}}$ are 
\[ p_{0,\; \varphi-1} \textrm{ to  } p_{{\varphi -1 \choose 2},\; \varphi -1},\textrm{ inclusive.} \]
\end{lemm}
\begin{pf} Observe that $\delta_G\geq 1$ for every $G\in \mathcal{G}$ (with the exception of $K_1$). Recall that a tree is a minimally connected graph and that a tree with $n$ vertices has exactly $n-1$ edges. Therefore, the minimum possible $\epsilon$ for a degree sequence $s$ of any connected graph is $\epsilon=\frac{2\varphi-2}{2}=\varphi-1$, examples of which include, $\underbrace{\{\varphi-1,1,\dots,1\}}_{\varphi}$, realised as a star graph, and $ \underbrace{\{2,\dots,2,1,1\}}_{\varphi}$, realised as a path graph. It follows from Lemma \ref{Ent} that the uppermost non-empty entry in column $j=\varphi -1$ of $P_{\mathcal{G}}$ is contained in row $\min\{\epsilon\} -\varphi +1=(\varphi -1)-\varphi +1=0$. \\

Observe that for any $G\in \mathcal{G}$ where $d(G)=s$ then $\Delta_G\leq \varphi-1$ as $G$ is simple. It follows that the maximum possible $\epsilon$ for any connected graph is $\epsilon={\varphi \choose 2}$ i.e. $s=\underbrace{\{\varphi-1,...,\varphi-1\}}_{\varphi} = d(K_{\varphi})$. It follows from Lemma \ref{Ent} that the lowermost non-empty entry in column $j=\varphi -1$ of $P_{\mathcal{G}}$ is contained in row $\max\{\epsilon\} -\varphi +1={\varphi \choose 2}-\varphi +1=\frac{\varphi^2-3\varphi +2}{2}={\varphi-1 \choose 2}$. This argument is summarised in Figure \ref{table}.\\

 \begin{figure}[h]
{\renewcommand{\arraystretch}{1.5} 
\[ \begin{array}{c|ll|c}
 p_{i,j}\textrm{ in } P_{\mathcal{G}} & \textrm{Sequences }\{s_1,\dots, s_n\} &  & \epsilon \\ \hline
p_{\;0,\;\varphi -1} & \{n-1,1,\dots,1\}, \hspace{0.5cm} \cdots  & , \{2,\dots,2,1,1\} & \varphi -1 \\
\vdots & \hspace{1cm}\vdots  &   & \vdots \\
p_{\;{\varphi-1 \choose 2},\; \varphi -1} & \{n-1,\dots,n-1\} & &  {\varphi \choose 2}  \\
\end{array}  \] }
\caption{All possible connected degree sequences of length $n$.}
\label{table}
\end{figure}

\end{pf}

\begin{crry}\label{EntCorr} The non-empty entries in column $j$ of $P_{\mathcal{G}}$ are 
\[p_{0,\; j} \textrm{ to } p_{\frac{j^2 - j}{2},\; j}, \textrm{ inclusive.}\]
\end{crry}
\begin{pf} The result follows from Lemma \ref{EmpEnt} by letting $\varphi=j+1$.
\end{pf}
It follows that column $j$ in $P_{\mathcal{G}}$ has $\frac{j^2 - j}{2}+1$ non-empty entries, for example, column $j=3$ has $\frac{3^2 - 3}{2}+1=4$ non-empty entries as shown in Figure \ref{Gamma}.

\section{Results}

\begin{thrm}\label{Main} Given a sequence $s=\{s_1,...,s_n\}$ of positive integers, with the associated pair $(\varphi, \epsilon)$, such that $s_i\geq s_{i+1}$ for $i=1,...,n-1$ then $s$ is connected if and only if 
\begin{itemize}
\item $\epsilon \in \mathbb{N}$,
\item $\varphi -1\leq \epsilon \leq {{\varphi \choose 2}}$, 
\item $s_1\leq \varphi -1$ and $s_n \geq 1$.
\end{itemize}
\end{thrm}

\begin{pf}  ($\Rightarrow$) Clearly $\epsilon \in \mathbb{N}$ is a necessary condition for any sequence $s$ to be realisable as half the sum of the degrees in any graph is the number of edges in that graph and this must be a natural number. The necessity of the condition $\varphi -1\leq \epsilon \leq {{\varphi \choose 2}}$ follows from Lemma \ref{EmpEnt}. The necessity of the condition $s_1\leq \varphi -1$ follows directly from the definition of a simple graph. Note that $s_n\geq 1$ as otherwise any $s_i=0$ (apart from the trivial case where $s=(0)$) would result in a graph containing at least one isolated vertex which would be disconnected.\\

($\Leftarrow$) Suppose that $s=\{s_1,...,s_n\}$ is connected. This means that $s$ is the degree sequence of a connected graph $G$, hence $\sum\limits_{i=1}^{n}deg(v_i)=2|V_G|$ and so $\epsilon\in\mathbb{N}$. A tree is by definition a minimal connected graph and every tree with $n$ vertices contains exactly $n-1$ edges, hence $\epsilon\geq \varphi -1$. As $G$ is simple then $deg(v_i)\leq n-1$ for all $i=1,...,n$. The maximal simple graph on $n$ vertices is the complete graph $K_n$, which has the degree sequence $\{n-1,...,n-1\}$ and $|E_{K_n}|={n \choose 2}$, hence $s_1\leq n-1$ and $\epsilon \leq {\varphi \choose 2}$.   
\end{pf}

Before stating the next result the following definition is required. 

\begin{dfn} A finite sequence $s=\{s_1,...,s_n\}$ of positive integers is called {\it necessarily connected} if $s$ can only be realisable as a connected simple graph.
\end{dfn}

\begin{thrm}\label{Crry1} Given a sequence $s=\{s_1,...,s_n\}$ of positive integers, with the associated pair $(\varphi, \epsilon)$, such that $s_i\geq s_{i+1}$ for $i=1,...,n-1$ then $s$ is necessarily connected if and only if $s$ is connected and 
$\epsilon > {{\varphi -2 \choose 2}} +1.$
\end{thrm}

\begin{pf}  ($\Rightarrow$) Clearly it is necessary for $s$ to be connected if it is to be necessarily connected. It is required to show that it is necessary for $ \epsilon > {{\varphi -2 \choose 2}} +1$. Consider a sequence $s=\{s_1,...,s_n\}$ such that $ \epsilon = {{\varphi -2 \choose 2}} +1.$ Observe that one such sequence is $s'=\{n-3,...,n-3,1,1\}$ which has $(\varphi(s'),\epsilon(s'))=\left(n, \frac{(n-2)(n-3)+2}{2}\right)=\left(n,{n -2 \choose 2} +1\right)$. Observe that $\{n-3,...,n-3,1,1\}=d(G_1)$, see Figure \ref{graphs}, where $G_1\simeq K_{n-2}\sqcup K_2$ with $V_{K_{n-2}}=\{v_1,...,v_{n-2}\}$ and $V_{K_2}=\{v_{n-1},v_n\}$. However, $s'=\{n-3,...,n-3,1,1\}$ is, in fact,  connected as $s'$ is also the degree sequence of the connected graph $G_2$ shown in Figure \ref{graphs}. Note that $v_iv_j\in E_{G_1}$ but $v_iv_j\not\in E_{G_2}$. Therefore, it is necessary that $ \epsilon > {{\varphi -2 \choose 2}} +1$ if $s$ is to be necessarily connected as the sequence $s'=\{n-3,...,n-3,1,1\}$, with $\epsilon(s')={{\varphi -2 \choose 2}} +1$, is realisable as a disconnected graph.

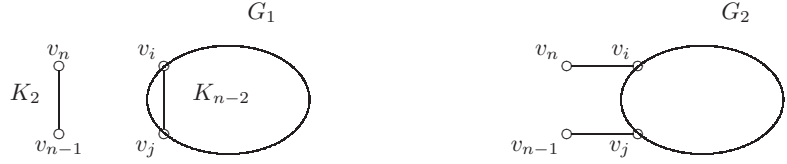
\begin{figure}[h]
\begin{center}
\scalebox{0.9}{$\begin{xy}\POS (-5,5) *\cir<2pt>{} ="a" *+!D{v_n},
(10.5,5) *\cir<2pt>{} ="b" *+!DR{v_i},
 (-5,-5) *\cir<2pt>{} ="c" *+!U{v_{n-1}},
(10.5,-5) *\cir<2pt>{} ="d" *+!UR{v_j},
  (25,12)*+!{G_1},
  (19,0)*+!{K_{n-2}},
   (-10,0)*+!{K_{2}},
  
\POS "a" \ar@{-}  "c",
\POS "b" \ar@{-}  "d",

\POS(20,0),  {\ellipse(12,8)<>{}},

\POS (70,5) *\cir<2pt>{} ="a" *+!DR{v_n},
(80.5,5) *\cir<2pt>{} ="b" *+!DR{v_i},
 (70,-5) *\cir<2pt>{} ="c" *+!UR{v_{n-1}},
(80.5,-5) *\cir<2pt>{} ="d" *+!UR{v_j},
  (95,12)*+!{G_2},
  
\POS "a" \ar@{-}  "b",
\POS "c" \ar@{-}  "d",

\POS(90,0),  {\ellipse(12,8)<>{}},

 \end{xy}$}

\caption{ $d(G_1)=d(G_2)=\{n-3,...,n-3,1,1\}$.}
\label{graphs}
\end{center}
\end{figure}

($\Leftarrow$) It is now required to show that if $s$ is connected and $ \epsilon > {{\varphi -2 \choose 2}} +1$ then $s$ is necessarily connected. To show this it is required to show that the maximum number of edges in a graph with $n$ vertices which can be disconnected is ${n-2 \choose 2} +1$. The graph $G_1$ in Figure \ref{graphs} shows that such a disconnected graph exists, so it is now required to show that a graph with $ \epsilon = {{\varphi -2 \choose 2}} +1$ is maximal i.e. adding one edge always results in a connected graph. \\

Observe that any maximal disconnected graph on $n$ vertices necessarily contains two disjoint connected components $H_1$ and $H_2$. To maximise the number of edges in $H_1 \sqcup H_2$ then $H_1\simeq K_a$ and $H_2\simeq K_b$ where $a+b=n$. So, maximising $|E_{H_1}|+|E_{H_2}|$ is equivalent to minimising the number of edges in a complete bipartite graph $K_{a,b}$ as $K_{n}\setminus (E_{H_1\sqcup H_2}) \simeq K_{a,b}$.\\

Let $a+b=n$, with $a\leq b$, then $|E_{K_{a,b}}|=ab$ where $a,b\in\{1,...,n-1\}$. Note that $a>0$ as $G$ is disconnected i.e. $K_a\neq K_0=(\varnothing, \varnothing)$. It is straightforward to show that $ab$ attains its maximum at $a=b=\frac{n}{2}$, when $n$ is even, and at $a=\lfloor\frac{n}{2}\rfloor, b=\lceil\frac{n}{2}\rceil$ when $n$ is odd. It follows that $ab$ is minimised when $a=1$ and $b=n-1$. However, observe that $a>1$ as $a=1$ implies that $H_1\simeq K_1$ which means that $d(G)$ contains a term equal to zero which contradicts the $s_n\geq 1$ condition. Hence $|E_{K_{a,b}}|$, with $a+b=n$, is minimised when $a=2$ and $b=n-2$ and so the maximal disconnected graph on $n$ vertices is isomorphic to $K_{2}\sqcup K_{n-2}$, see $G_1$ in Figure \ref{graphs}. Notice that the union of $G_1$ and any edge in $\overline{G_1}$, the complement of $G_1$, results in a connected graph. 
\end{pf}

\begin{crry} All simple graphs with $n$ vertices and at least $\frac{n^2-5n +10}{2}$ edges are connected.
\end{crry}
\begin{pf}  As shown in Theorem \ref{Crry1}, a maximal disconnected graph with $n$ vertices is isomorphic to $K_{n-2}\sqcup K_2$ and thus has $\left({n -2 \choose 2} +1\right)$ edges. It follows that any simple graph with $n$ vertices and at least ${{n -2 \choose 2}} +2 = \frac{(n-2)(n-3)+4}{2}=\frac{n^2-5n+10}{2}$ edges is connected. 
\end{pf}
Note that for all $n\in \mathbb{N}$, $n^2-5n$ is even and $n^2-5n +10>0$, hence $\frac{n^2-5n +10}{2}\in \mathbb{N}$. 



\bibliographystyle{plain}      
\bibliography{refs2015}

\end{document}